\documentclass[smallextended]{svjour3} 

\RequirePackage{fix-cm}
\smartqed  
\usepackage{graphicx, amsmath, adjustbox, multirow}

\journalname{Biomechanics and Modeling in Mechanobiology}
\usepackage{color}

\usepackage{bm}
\usepackage{url}

\begin{document}

\title{A computational study of aortic reconstruction in single ventricle patients}

%

\titlerunning{Aortic reconstruction in single ventricle patients}

\author{Alyssa M. Taylor-LaPole \and Mitchel J. Colebank \and Justin D. Weigand \and Mette S. Olufsen \and Charles Puelz}


\institute{Alyssa M. Taylor-LaPole \at
              Department of Mathematics\\
              North Carolina State University\\
              Raleigh, NC 27695\\
              \email amtayl25@ncsu.edu   
           \and
           Mitchel J. Colebank \at
              Edwards Lifesciences Foundation Cardiovascular Innovation and Research Center \\
              Department of Biomedical Engineering\\
              University of California, Irvine\\
              Irvine, CA 92617\\
              \email mcoleban@uci.edu
           \and
           Justin D. Weigand \at
              Department of Pediatrics, Division of Cardiology\\
              Baylor College of Medicine and Texas Children's Hospital\\
              Houston, TX 77030\\
              \email Justin.Weigand@bcm.edu
           \and 
           Mette S. Olufsen \at
              Department of Mathematics \\
              North Carolina State University\\
              Raleigh, NC 27695\\
              \email msolufse@ncsu.edu 
            \and
           Charles Puelz \at
              Department of Pediatrics, Division of Cardiology\\
              Baylor College of Medicine and Texas Children's Hospital\\
              Houston, TX 77030\\
              \email charles.puelz@bcm.edu}

\date{Received: date / Accepted: date}

\maketitle

\begin{abstract}
Patients with hypoplastic left heart syndrome (HLHS) are born with an underdeveloped left heart. They typically receive a sequence of surgeries that result in a single ventricle physiology called the Fontan circulation. While these patients usually survive into early adulthood, they are at risk for medical complications, partially due to their lower than normal cardiac output, which leads to insufficient cerebral and gut perfusion. While clinical imaging data can provide detailed insight into cardiovascular function within the imaged region, it is difficult to use this data for assessing deficiencies in the rest of the body and for deriving blood pressure dynamics.  Data from patients used in this paper include three dimensional, magnetic resonance angiograms (MRA), time-resolved phase contrast cardiac magnetic resonance images (4D-MRI) and sphygmomanometer blood pressure measurements. The 4D-MRI images provide detailed insight into velocity and flow in vessels within the imaged region, but they cannot predict flow in the rest of the body, nor do they provide values of blood pressure.  To remedy these limitations, this study combines the MRA, 4D-MRI, and pressure data  with 1D fluid dynamics models to predict hemodynamics in the major systemic arteries, including the cerebral and gut vasculature. A specific focus is placed on studying the impact of aortic reconstruction occurring during the first surgery that results in abnormal vessel morphology. To study these effects, we compare simulations for an HLHS patient with simulations for a matched control patient that has double outlet right ventricle (DORV) physiology with a native aorta. Our results show that the HLHS patient has hypertensive pressures in the brain as well as reduced flow to the gut. Wave-intensity analysis suggests that the HLHS patient has irregular circulatory function during light upright exercise conditions and that predicted wall-shear stresses are lower than normal, suggesting the HLHS patient may have hypertension.

\keywords{hypoplastic left heart syndrome \and Fontan circulation \and double outlet right ventricle \and perfusion \and fluid mechanics \and simulation \and systemic circulation}
\end{abstract}

\section{Introduction} 

Hypoplastic left heart syndrome (HLHS) is a congenital disease characterized by an underdeveloped left heart. For most of these patients, the left ventricle is non-functional \cite{Tworetzky2001}, i.e.~it is unable to generate sufficient cardiac output to perfuse the systemic vasculature with oxygen-rich blood. If left untreated, infants with this defect do not survive. The optimal treatment is a heart transplant, but since few infant hearts are available, patients are usually treated via a sequence of palliative surgeries that result in a single ventricle physiology called the Fontan circulation. However, most HLHS patients will eventually require a heart transplant due to failure of the Fontan circuit \cite{Deal2012}. These patients live with a single ventricle pump and therefore do not have the typical capacity for transporting blood throughout the body. Moreover, surgically reconstructed vessels degenerate over time, further reducing the ability of the single ventricle to generate sufficient power \cite{Mahle1998,Voges2010}. The sequence of surgical procedures recommended for this patient population are an attempt to reconstruct and rearrange vessels in the vicinity of the heart so that the single ventricle is functional. These surgeries are most commonly called the Norwood, Glenn, and Fontan procedures \cite{Voges2010}.

The {\em Norwood procedure} is performed during the first few weeks of life. This procedure involves the construction of a new aorta, herein referred to as the {\em reconstructed aorta}. This vessel connects the main pulmonary artery and diminutive native aorta to create a connection between the single functioning ventricle and systemic circulation. The surgical reconstruction involves the addition of a homograph patch comprised of tissue harvested from the main pulmonary artery \cite{Mahle1998}. Since the reconstructed aorta is attached to the main pulmonary artery, it is necessary to add a shunt to transport blood to the pulmonary circulation for reoxygenation. Two shunts are commonly used: a Blalock-Taussig-Thomas shunt, which connects the subclavian artery or carotid artery to the pulmonary artery, or a Sano shunt, which connects the single ventricle directly to the pulmonary vasculature \cite{Ohye2010}.

The {\em Glenn procedure} is typically performed six months after birth. In this surgery, the superior vena cava is connected to the main pulmonary artery in order to transport venous blood returning from the upper body to the lungs. This new pathway eliminates the need for the shunt placed during the Norwood procedure and it is therefore removed. The resulting circulation mixes oxygenated blood from the pulmonary circulation with deoxygenated blood from the lower body. This mixed blood is subsequently transported to the systemic arteries \cite{Gobergs2016}. 

In the {\em Fontan procedure}, performed at age 18-36 months, the inferior vena cava is connected to the main pulmonary artery, allowing blood from the lower body to also travel to the lungs. The final physiology is characterized by serialized pulmonary and systemic circulations supported by the single functioning ventricle \cite{Fontan1971}.  

Many HLHS patients live with a single ventricle circulation until early adulthood, but they often experience serious complications, partly as a result of reduced cardiac output \cite{Gewillig2016}. This leads to insufficient cerebral and gut perfusion comorbid with increased risk of stroke \cite{Saiki2014} and chronic venous congestion contributing to Fontan-associated liver disease (FALD) \cite{Gordon2019,Navaratnam2016}. We hypothesize two important causes of complications: (1) impaired hemodynamic transport driven by the single ventricle and (2) degeneration in the reconstructed vessels over time \cite{Biglino2012,Gewillig2016,Mahle1998,Mitchell2018,Revell2013}. We hypothesize that the vessel experiencing the most significant degeneration is the reconstructed aorta. Examination of images (used in this study) clearly shows the aorta in HLHS patients is significantly larger than for the control group comprised of DORV patients. An increase in the aortic cross-sectional area will decrease pressure, requiring the functioning ventricle to pump harder in order to perfuse all essential organs. A cause for aortic degeneration could be the noncompliant patch material harvested from the main pulmonary artery, which does not continue to grow along with the native aortic tissue \cite{Mahle1998}. This might result in vessel remodeling, leading to an increase in vessel stiffness \cite{Ou2008} and cross-sectional area \cite{Voges2015}. These changes would likely lead to increased pulse pressure and the formation of vortices, particularly in the diastolic phase of the cardiac cycle, where flow reversal has been observed \cite{Mitchell2004,Mitchell2018}. 

To examine the impact of aortic remodeling in HLHS patients, we use computational simulations in order to compare hemodynamic predictions in an HLHS patient and double outlet right ventricle (DORV) patient used as the control. The latter group has single ventricle physiology but a native aorta, so they provide an ideal control group. Treatments for single ventricle physiology have been investigated heavily in experimental studies. Cardis et al. \cite{Cardis2006} examined the impact of aortic reconstruction on elastance properties of the aorta in HLHS, DORV, and other single ventricle patients. Their results demonstrated that reconstructed aortas in HLHS patients had lower distensibility and higher stiffness than the other single ventricle patients. Other experimental studies have examined the effectiveness and timing of the Fontan surgeries \cite{Mainwaring1994,Tanoue2001,Yagi2017}. For example, the study by Tanoue et al. \cite{Tanoue2001} compared results from 18 patients who had the Glenn procedure before the Fontan procedure, with 23 patients who proceeded directly to the total cardiopulmonary connection. The authors found that conducting the Glenn procedure first (which is currently recommended) significantly improved patient outcomes. Similar results were found in the studies by Mainwaring et al. \cite{Mainwaring1994}, who examined concentrations of hormones associated with fluid retention, and Yagi et al. \cite{Yagi2017}, who measured cerebral oxygenation. Other studies focused on failure of the Fontan circulation. For example, Kotan et al. \cite{Kotan2018} studied health records from 500 Fontan patients and found that premature morbidity and death were caused by circulatory failure, multiorgan failure, cerebral/renal issues, and pulmonary failure. They also found that lifespan increased with early intervention and treatment.  Saiki et al. \cite{Saiki2014} performed a statistical analysis of wave intensity in the carotid artery in 34 patients with a Fontan circulation and 20 controls. Their results showed that carotid blood flow was much lower in reconstructed Fontan patients than in healthy, double ventricle controls, suggesting that Fontan patients have lower cerebral perfusion.  


Most computational studies examining the Fontan circulation (and its preceding physiologies) use three-dimensional (3D) computational fluid dynamics (CFD) models in order to predict velocity distributions and power losses, mainly through the surgically created connection between the vena cavae and pulmonary arteries \cite{Marsden07,Marsden09,Ahmed2021}. Such simulations are ideal for assessing complex velocity patterns, and have been used to explore the effectiveness of various types of Fontan procedures \cite{Bove2007,Pekkan2009,Prather2022,Sun2014}.  Bove et al. \cite{Bove2007} and Pekkan et al. \cite{Pekkan2009} focused on predicting flow distribution and power loss by comparing two second-stage surgeries, the hemi-Fontan and Glenn Procedures. Pekken et al. \cite{Pekkan2009} used idealized and patient-specific models and emphasized the importance of examining patient-specific geometries, since results varied significantly between patients, particularly in reconstructed vessels with sharp turns (since these induce high resistance to flow).  Recently, Prather et al. \cite{Prather2022} found that in Fontan circuits at risk of failure, the addition of an ``injection jet" shunt (IJS) that draws blood from the aortic arch to lower inferior vena cava can lower venous pressure. Other important computational work includes that of Bazilevs et al. which compares CFD and fluid-structure interaction models of the total cavopulmonary connection in the Fontan circulation \cite{Bazilevs2009}. The authors concluded that the additional detail provided by the fluid-structure interaction model was important for describing hemodynamics in this anatomy.

An alternative computational approach involves one-dimensional (1D) fluid dynamics models, which can efficiently predict wave propagation in large networks as well as account for fluid-solid interaction. For example, Puelz et al. \cite{Puelz2017} used 1D models to investigate differences between a Fontan circulation with either a fenestration or a hepatic vein exclusion. Their study combined 1D arterial and venous networks with zero-dimensional (0D) heart and organ bed models. They found that both modifications to the Fontan circuit improve flow to the gut. 

\begin{figure}[b!]
    \centering
    \includegraphics[width=0.75\textwidth]{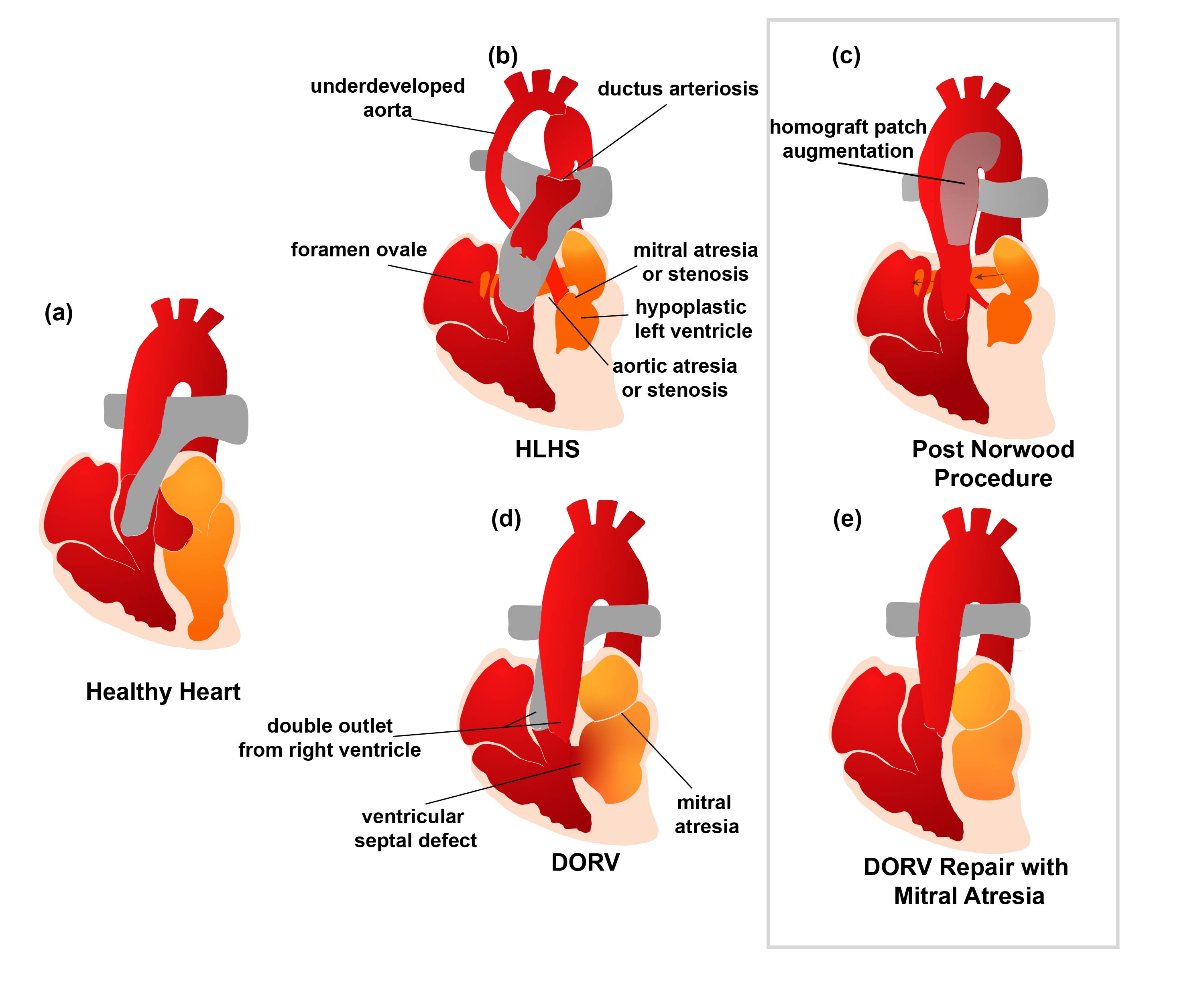}
    \caption{A schematic comparison of the physiology of a healthy heart, an HLHS heart, a DORV heart, and the HLHS and DORV hearts after the Norwood procedure. \textbf{(a)} A healthy heart with the aorta attached to the left ventricle and the pulmonary artery to the right ventricle. \textbf{(b)} The HLHS heart at birth has an underdeveloped aorta, an underdeveloped left ventricle, an extra vessel connecting the pulmonary artery, and the aorta (the ductus arteriosis), and a hole connecting the left and right atria (the foramen ovale). \textbf{(c)} During the Norwood procedure, a homograft patch augmentation is used to construct a new aorta (termed the reconstructed aorta), which is then attached to the functioning single ventricle. The native underdeveloped aorta emanating from the left ventricle is kept in order to perfuse the coronary arteries. \textbf{(d)} The DORV heart at birth has both the pulmonary artery and aorta attached to the right ventricle. \textbf{(e)} DORV physiology after the first procedure. No aortic intervention is necessary for DORV single ventricle variants since the native aorta suffices to provide systemic output. The gray box highlights the major difference (the size of the aorta) between the two patients studied.}
    \label{fig:HLHS_DORV}
\end{figure}

The experimental and computational studies discussed above provide significant insight into complex single ventricle phyiology, demonstrating their importance in enhancing our understanding of treatments for HLHS. To our knowledge, there are no computational studies on the effect of aortic reconstruction in HLHS patients that incorporate a comparison to DORV control patients. The goal of this study is to create 1D patient-specific vessel network models that predict blood pressures and flows in regions of interest as well as perfusion to essential organs.  Patient-specific networks are constructed from MRA images of two age- and size-matched children; one is an HLHS patient with a reconstructed aorta and the other is a DORV patient with a native aorta. Both patients have a single ventricle Fontan circulation consisting of the systemic organs and lungs in series (Figure \ref{fig:HLHS_DORV}). We construct model geometries from the MRA data, calibrate models to the 4D-MRI and blood pressure data, and then use the calibrated models to predict the impact of aortic reconstruction on cerebral and gut perfusion, blood pressure, wave intensity, and shear stress.

\section{Methods}

Vessel dimensions within the imaged region for the HLHS and DORV patients are extracted from magnetic resonance angiography (MRA) images \cite{Stankovic2014}. To predict the desired quantities, these networks are extended to include all major systemic arteries in the body and brain. To obtain a model that fits to data, vessels are scaled to patient weight. A one-dimensional (1D) fluid dynamics model is solved in these vessel networks to predict blood flow, pressure, wave intensity, and vessel wall shear stress.

\subsection{Data and Network Geometry}

This study uses retrospective hemodynamic measurements from one HLHS and one age- and size-matched DORV control patient. The HLHS patient received a Norwood reconstructed aorta, and both patients had Glenn and Fontan surgeries at the Texas Children's Hospital Heart Center in Houston, TX. Data collection for this study was approved by the Baylor College of Medicine Institutional Review Board (H-46224: ``Four-Dimensional Flow Cardiovascular Magnetic Resonance for the Assessment of Aortic Arch Properties in Single Ventricle Patients''). Patient characteristics, listed in Table \ref{table:Patient}, include age, height, gender, weight, average resting blood pressure (systolic and diastolic), and cardiac output. 

\begin{table}[b!]
\centering
\caption{Age, height, weight, blood pressures, and cardiac cycle duration for the DORV and HLHS patients. Cuff pressures were obtained using a sphygmomanometer from the upper arm in the supine position. These pressures are compared to systolic and diastolic pressures of the brachial artery in the supine rest position. The cardiac cycle duration was measured at the time of the 4D-MRI acquisition in the supine rest position.}
\label{table:Patient}       
\begin{tabular}{l|rr}
\hline\noalign{\smallskip}
 & \textbf{DORV} & \textbf{HLHS}  \\
\noalign{\smallskip}\hline\noalign{\smallskip}
\textbf{Age (yrs)} & 12 & 11 \\
\textbf{Height (cm)} & 154.3 & 151.4 \\
\textbf{Weight (kg)} & 59.6 & 62.0 \\
\textbf{Systolic pressure (mmHg)} & 110 & 116 \\
\textbf{Diastolic dressure (mmHg)} & 67 & 65 \\
\textbf{Cardiac cycle length (s)}& 0.658 & 0.615\\
\noalign{\smallskip}\hline
\end{tabular}
\end{table}

\subsubsection{Measurements} \label{Sec:Data}

Patient data are extracted from high spatial resolution three-dimensional (3D) contrast enhanced magnetic resonance angiography (MRA) images, three-dimensional, time-resolved phase contrast cardiac magnetic resonance (4D-MRI) images, and sphygmomanometer blood pressure measurements. Imaging data include the ascending aorta, aortic arch, brachiocephalic vessels, and descending thoracic aorta. Imaging studies were performed using a 1.5T Siemens Aera magnet (Siemens Healthineers, Erlangen, Germany). 

Using localizing sequences from the cardiac MRI, time-resolved, contrast enhanced dynamic MRA was performed with bolus injection of 0.1 mL/kg of intravenous gadolinium contrast. MRA was performed in the sagittal plane over a 30-60 second breath-hold, with 3-7 measurements acquired after contrast administration at a temporal resolution of 3-5 seconds/measurement based on heart rate.  90 to 120 slices per measurement were acquired with a slice thickness of 1.2 to 1.4 mm and reconstructed voxel dimensions of 1.2$\times$1.2$\times$1.2 mm. Images were stored on a password protected server in DICOM format.

\begin{figure}[t!]
    \centering
    \includegraphics[width=0.75\textwidth]{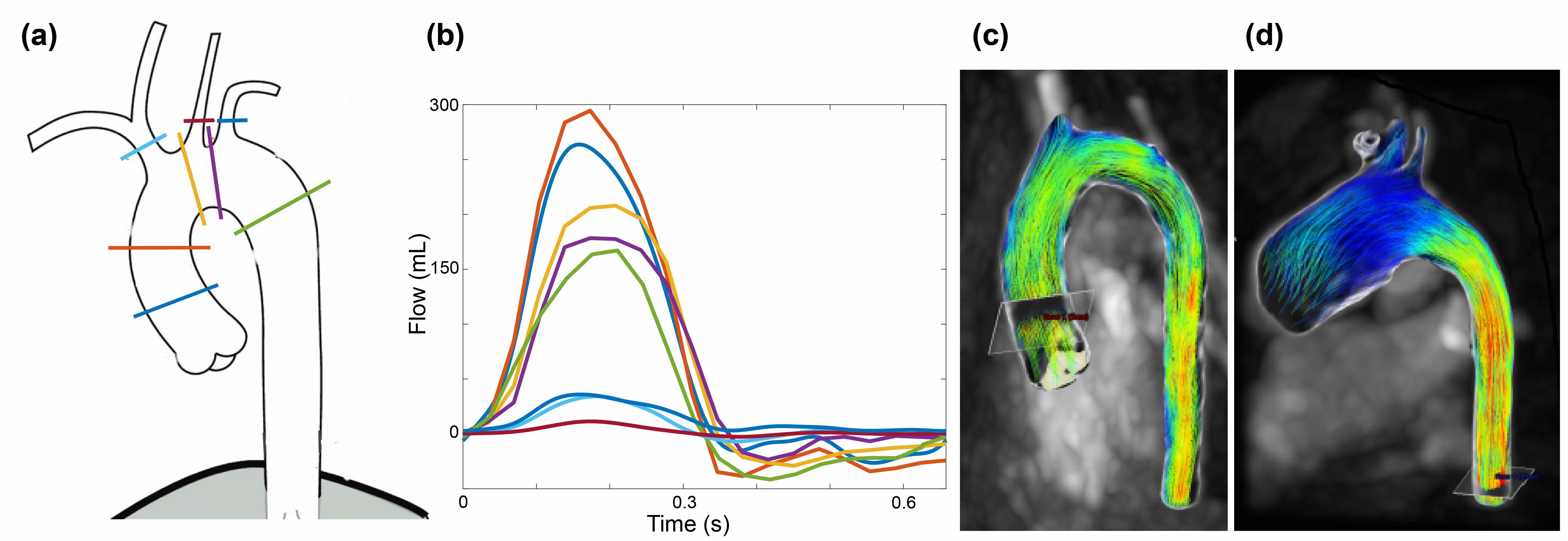}
    \caption{\textbf{(a)} A sketch of the aorta and the different planes on which the flow data was measured. \textbf{(b)} Flow profiles plotted for the DORV patient corresponding to the planes along the vessels in \textbf{(a)}. \textbf{(c)} Flow velocity streamlines in the DORV patient. \textbf{(d)} Flow velocity streamlines in the HLHS patient.}
    \label{fig:PatientProtocol}
\end{figure}

Three-dimensional time resolved phase contrast CMR (4D-MRI) was performed using the gradient echo phase contrast sequence that utilized both ECG and respiratory navigation to allow for free-breathing acquisition. The sequence was prescribed as a sagittal acquisition to cover the entire thoracic aorta and proximal brachiocephalic vessels. The 4D flow sequence was acquired with the patient free breathing with a slice thickness of 1 to 2.5 mm. Velocity encoding was set to 10\% above the highest velocity expected in the aorta. The accuracy of the 4D-MRI sequences are in part determined by spatial and temporal resolutions as well as an adequate signal-to-noise ratio. In an effort to maximize the latter, the 4D-MRI sequence was performed after the acquisition of the MRA since the contrast provided an enhanced signal-to-noise ratio for the 4D-MRI sequence.

4D-MRI post processing was performed using CVI 42 software (Circle Cardiovascular Imaging, Calgary, CA). The thoracic aorta and proximal brachiocephalic vessels were included in the region of interest. A velocity mask of the aorta was generated and assessed for aliasing in three orthogonal planes. No aliasing was noted on the studies and the aorta was segmented from the surrounding structures. A centerline was established through the aorta and subsequently in the brachiocephalic, subclavian, and carotid arteries. Flow plane slices at specified locations along the aorta were constructed orthogonally to the centerline. For each prescribed flow plane location, a flow waveform was generated by integrating the masked velocity field over the plane. For this study, we extracted 9 volumetric flow waveforms along the ascending aorta, aortic arch, descending aorta, thoracic aorta, brachiocephalic, left common carotid, and left subclavian arteries (refer to Figure \ref{fig:PatientProtocol}).

\begin{figure}[t!]
    \centering
    \includegraphics[width=1.0\textwidth]{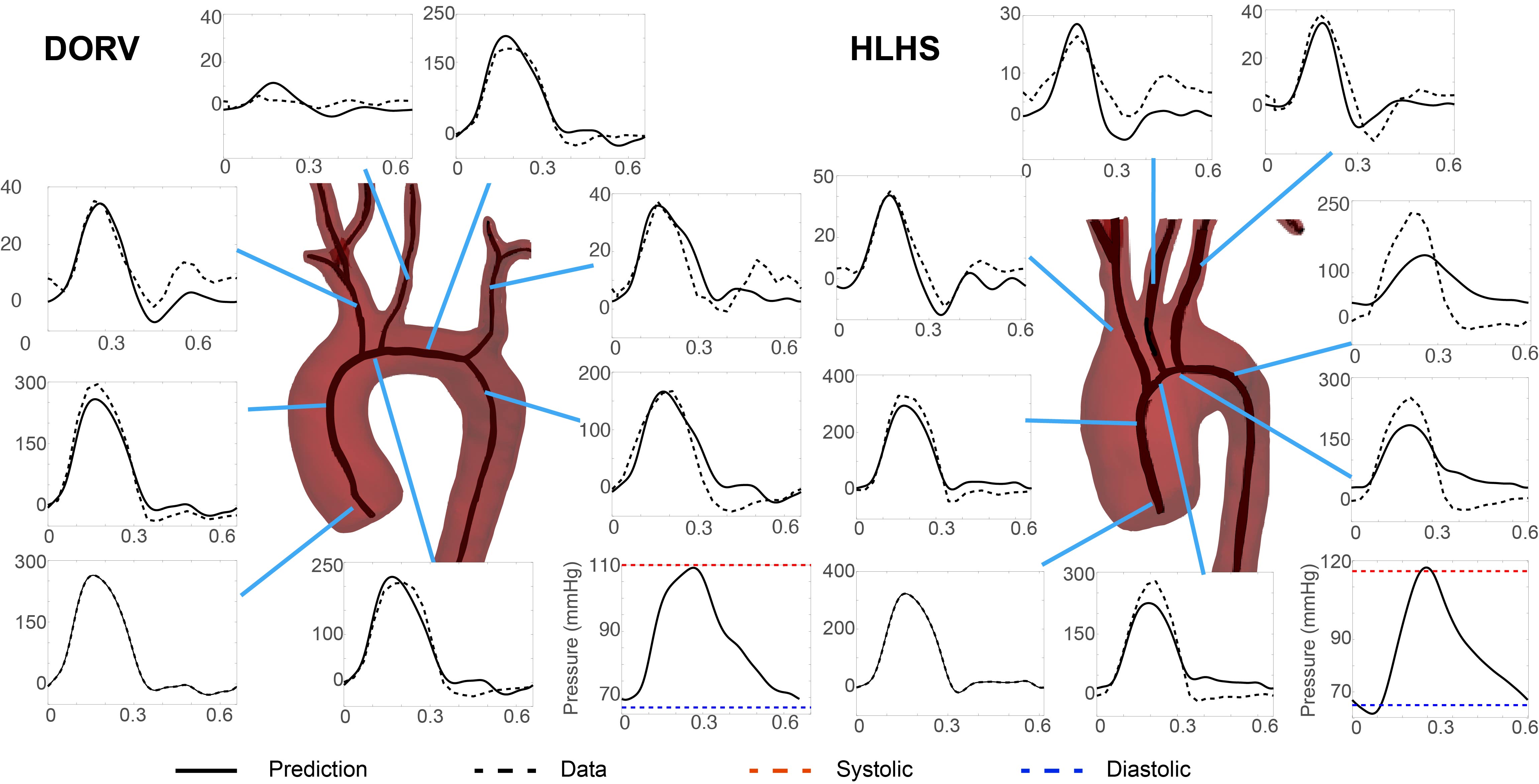}
    \caption{A visualization of the flow and pressure waveforms. Flows derived from the 4D-MRI data are shown as dashed lines and flows predicted by the models are shown as solid lines. Measured flow data are scaled to ensure conservation of flow required for the 1D model predictions. Cuff pressure measurements are obtained for each patient and shown in the bottom right panels in red (systolic) and blue (diastolic). These pressure measurements are taken in the upper left arm and the model pressure predictions, from the upper brachial artery, are shown in the black solid line.}
    \label{fig:flowdata}
\end{figure}

Systolic and diastolic cuff pressures and resting heart rates (which were used to determine the length of one cardiac cycle) were also measured. To ensure conservation of blood flow, the brachiocephalic, subclavian, and left common carotid flow waveforms were scaled. We chose to scale these flow waveforms since the corresponding vessels are smaller (compared to the aorta), and a small uncertainty in the cross-sectional area defined by the velocity mask gives rise to higher uncertainty in the flow predictions. In addition, we scaled the flow measurements in the aorta if the average flow was bigger than the flow in the preceding slice. Average flow values before and after scaling are given in Table \ref{tab:conservation} and the scaled flow waveforms are shown as dashed lines in Figure \ref{fig:flowdata}. Flows were scaled by multiplying the original waveform by a certain percentage in order to enforce mass conservation. Average flow values were calculated by integrating waveforms over the duration of the patient's cardiac cycle and dividing by that duration.

\begin{table}[b]
\caption{Average flow data from the 4D-MRI images before and after scaling to ensure mass conservation. Values listed are in L/min.}
\centering
\label{tab:conservation}   
\begin{tabular}{l|ll|ll}
\hline\noalign{\smallskip}
& \multicolumn{2}{|c|}{\textbf{DORV}} & \multicolumn{2}{|c}{\textbf{HLHS}}\\
\noalign{\smallskip}\hline\noalign{\smallskip}
 \textbf{Vessel} & \textbf{Data} & \textbf{Scaled} & \textbf{Data} & \textbf{Scaled}\\
\noalign{\smallskip}\hline\noalign{\smallskip}
Inflow & 4.06 & 4.06 & 5.08 & 5.08\\
Asc. Aorta & 4.14 & 4.06 & 4.49 & 5.08\\
Aortic Arch I & 3.69 & 3.32 & 4.75 & 4.51\\
Aortic Arch II & 4.22 & 3.13 & 4.16 & 4.12 \\
Thoracic Aorta & 2.95 & 2.62 & 3.67 & 3.67 \\
Brachiocephalic & 1.95 & 0.75 & 1.04 & 0.56 \\
L Comm. Carotid & 1.19 & 0.19 & 0.54 & 0.40 \\
L Subclavian & 0.82 & 0.50 & 0.74 & 0.45 \\
\noalign{\smallskip}\hline
\end{tabular}
\end{table}

\subsubsection{Geometric Domain}
\begin{figure}[b!]
    \centering
    \includegraphics[width=1.0\textwidth]{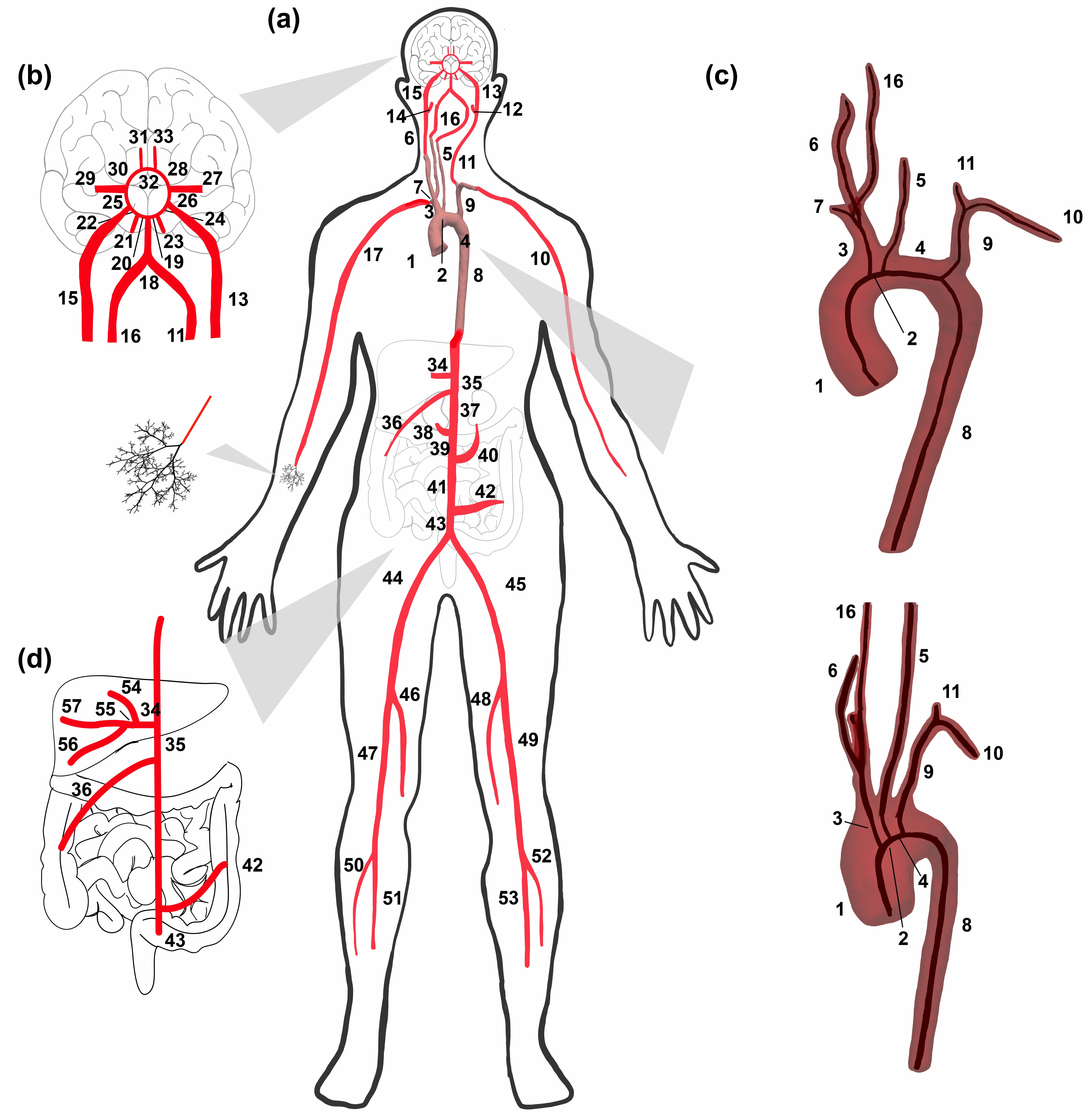}
    \caption{The network used for our patient specific model. The 3D rendered images of the aorta are segmented from the MRA data. \textbf{(a)} shows the entire network, including vessels outside the imaged region. The latter are calibrated using allometric scaling of literature values and patients weight. Attached to each terminal vessel is a structured tree model, shown on the right arm. \textbf{(b)} gives a detailed view of the cerebral circulation. \textbf{(c)} shows 3D renderings of the DORV and HLHS reconstructed patient aortas and nearby vasculature. Vessel centerlines are represented by the darker line through each segmentation. The upper image shows the DORV patient's aorta and the lower image shows the HLHS patient's reconstructed aorta. The ascending aorta from the HLHS patient has a much greater radius than that of the DORV patient's aorta. This difference is also apparent in Table \ref{table:Dim} with reference to the vessel radii.\textbf{(d)} shows a detailed view of the gut circulation included in this model.}
    \label{fig:Network}
\end{figure}

\begin{table}[t]
\caption{Vessel names and dimensions used in this study. Vessel numbers correspond to numbers given in Figure \ref{fig:Network}.}
\label{table:Dim}  
\adjustbox{max width = \textwidth}{
\begin{tabular}{ll|lll|lll}
\hline\noalign{\smallskip}
& & \multicolumn{3}{|c|}{\textbf{DORV}} & \multicolumn{3}{|c}{\textbf{HLHS}}\\
\noalign{\smallskip}\hline\noalign{\smallskip}
 \textbf{Number} & \textbf{Name} & $L$ \textbf{(cm)} & $R\text{in}$ \textbf{(cm)} & $R_\text{out}$ \textbf{(cm)} & $L$ \textbf{(cm)} & $R_\text{in}$ \textbf{(cm)} & $R_\text{out}$ \textbf{(cm)} \\
\noalign{\smallskip}\hline\noalign{\smallskip}
1 & Ascending aorta & 4.07 & 1.20 & 1.10 & 3.87 & 1.96 & 1.88 \\
2 & Aortic arch I & 1.95 & 1.10 & 1.10 & 1.93 & 1.60 & 1.20\\
3 & Brachiocephallic & 1.23 & 0.57 & 0.49 & 1.60 & 0.57 & 0.52 \\
4 & Aortic arch II & 1.94 & 0.95 & 0.88 & 3.77 & 1.55 & 1.06 \\
5 & L common carotid & 20.23 & 0.36 & 0.36 & 20.1 & 0.36 & 0.36\\
6,11 & Vertebral & 14.4 & 0.20 & 0.19 & 14.3 & 0.20 & 0.19\\
7,9 & Subclavian & 3.67 & 0.59 & 0.37 & 3.27 & 0.46 & 0.39\\
8 & Thoracic aorta & 15.17 & 0.88 & 0.70 & 15.07 & 1.44 & 0.69 \\
10,17 & Brachial & 20.23 & 0.35 & 0.30 & 20.1 & 0.34 & 0.30\\
12,14 & External carotid & 17.22 & 0.32 & 0.32 & 17.01 & 0.31 & 0.31\\
13,15 & Internal carotid I & 17.12 & 0.32 & 0.32 & 17.01 & 0.31 & 0.31\\
16 & R common carotid & 17.22 & 0.39 & 0.39 & 17.10 & 0.39 & 0.39\\
18 & Basilar & 2.76 & 0.15 & 0.15 & 2.74 & 0.15 & 0.15\\
19,20 & PCA I & 0.48 & 0.10 & 0.10 & 0.47 & 0.10 & 0.10\\
21,23 & PCA II & 8.18 & 0.10 & 0.10 & 8.13 & 0.10 & 0.10\\
22,24 & PCoA & 1.43 & 0.07 & 0.07 & 1.42 & 0.07 & 0.07\\
25,26 & Internal carotid II & 0.48 & 0.19 & 0.19 & 0.47 & 0.19 & 0.19\\
27,29 & MCA & 11.32 & 0.14 & 0.14 & 11.24 & 0.14 & 0.14\\
28,30 & ACA I & 1.14 & 0.11 & 0.11 & 1.13 & 0.11 & 0.11\\
31,33 & ACA II & 9.8 & 0.11 & 0.11 & 9.73 & 0.11 & 0.11\\
32 & ACoA & 0.29 & 0.07 & 0.07 & 0.28 & 0.07 & 0.07\\
34 & Celiac axis & 1.95 & 0.37 & 0.37 & 1.93 & 0.37 & 0.37\\
35 & Abdominal aorta & 5.16 & 0.59 & 0.57 & 5.12 & 0.59 & 0.47\\
36 & Superior mesenteric & 5.74 & 0.29 & 0.29 & 5.70 & 0.29 & 0.29\\
37 & Abdominal aorta & 0.97 & 0.57 & 0.55 & 0.97 & 0.57 & 0.55 \\
38,40 & Renal & 3.11 & 0.25 & 0.25 & 3.09 & 0.25 & 0.25\\
39 & Abdominal aorta & 0.97 & 0.55 & 0.53 & 0.97 & 0.55 & 0.53 \\
41 & Abdominal aorta & 10.31 & 0.53 & 0.51 & 10.24 & 0.53 & 0.50\\
42 & Inferior mesenteric & 4.86 & 0.16 & 0.16 & 4.83 & 0.15 & 0.15\\
43 & Abdominal aorta & 0.97 & 0.51 & 0.48 & 0.97 & 0.50 & 0.48 \\
44,45 & External iliac & 14.01 & 0.27 & 0.26 & 13.91 & 0.27 & 0.26\\
46,48 & Internal iliac & 4.86 & 0.26 & 0.26 & 4.83 & 0.26 & 0.26\\
47,49 & Femoral & 13.09 & 0.24 & 0.21 & 13.00 & 0.21 & 0.21\\
50,52 & Femoral & 43.09 & 0.21 & 0.21 & 42.81 & 0.21 & 0.21\\
51,53 & Deep femoral & 12.26 & 0.15 & 0.12 & 12.17 & 0.14 & 0.12\\
54 & Splenic & 5.99 & 0.21 & 0.19 & 5.95 & 0.21 & 0.19  \\
55 & Celiac axis & 1.90 & 0.25 & 0.25 & 1.89 & 0.25 & 0.25 \\
56 & Left Gastric & 6.75 & 0.15 & 0.14 & 6.71 & 0.15 & 0.14 \\
57 & Hepatic & 6.28 & 0.26 & 0.21 & 6.24 & 0.26 & 0.21 \\
\noalign{\smallskip}\hline
\end{tabular}}
\end{table}

For each patient, vessel dimensions are extracted from the MRA data. Images from the DORV and HLHS patients include the ascending aorta, aortic arch, brachiochephalic trunk, thoracic aorta, subclavian artery, common carotid artery, vertebral artery, and brachial artery. Vessel network extraction requires several steps. The first step is to generate a 3D volume rendering of the large vessels from the MRA image. The second step requires the generation of vessel centerlines, radii, and lengths. The third step involves the organization of these parameters  into a labeled graph that consists of edges corresponding to vessels, nodes corresponding to junctions, and radii and lengths for each vessel. Finally, the graph is augmented by attaching the main systemic arterial networks appearing outside the imaged region. The last step is done by scaling literature values of peripheral vessel networks to patient weight.

\paragraph{3D rendering.} Segmentation of the vessel geometries in the imaged region is performed using the open-source software 3DSlicer from Kitware Inc. \cite{Fedorov2012,Kikinis2014}. A 3D rendered volume is obtained using the built in thresholding\footnote{Image intensities are taken in the range 100 to 264 Houndsfield Units.}, cutting, and islanding  tools. The 3D geometry is saved to STL format and imported into Paraview (Kitware Inc. \cite{Utkarsh2015}), where it is converted to a VTK polygonal data file before processing using the Vascular Modeling Toolkit (VMTK, \url{http://www.vmtk.org/}) \cite{Antiga2008}. 

\paragraph{Centerlines.} VMTK generates centerlines from the 3D rendered volume by  inscribing spheres within each vessel \cite{Antiga2008}. At each point along the vessel, the radius is determined from the maximally inscribed sphere. VMTK generates centerlines that begin at a manually selected source point (the inlet of the aorta) and terminate at manually selected end points. In each vessel, the software then backtraces centerlines from the end points to the source point. Vessel junctions are set at the last point where two centerlines interact and are subsequently adjusted manually, since some junction points are placed before the barometric center of the junction. Additional details regarding methodology can be found in Antiga et al. \cite{Antiga2008}. Centerlines are saved as CSV files and imported into custom Matlab (The MathWorks, Inc., Natick, Massachusetts) software for postprocessing. Centerlines are shown in the 3D rendering of the aorta geometries in Figure \ref{fig:Network}, panel (c). The DORV and HLHS geometries are the top and bottom aortas respectively.

\paragraph{Directed graph.}
Previously developed algorithms implemented in Matlab are used to construct a labeled graph from the VMTK centerlines \cite{Colebank2019,Colebank2021}. Centerlines are converted to edges, corresponding to vessels, and nodes, corresponding to junctions. Each vessel is added to a connectivity matrix that determines the vessel network topology. Edges include $xyz$-coordinates along the centerlines and the corresponding radius determined from the maximally inscribed sphere within the vessel. Vessel length is calculated as a sum of the distances between $xyz$-coordinates. Vessel junctions are defined as the intersection of two centerlines, and terminal vessels (corresponding to sites for the structured tree models) are identified as branches with no distal daughter vessels. 

\paragraph{Whole body network.} Arteries outside of the imaged region are allometrically scaled from literature data according to patient weight. To ensure the beginning of the vessels equals the outlet radii of the patient specific vessel geometries created from the MRA data, we incorporate tapering of the large vessels as discussed below. Similar to previous studies by Melis et al. \cite{Melis2019} and Puelz et al. \cite{Puelz2017}, we used a 57-vessel network that includes all vessels captured in the MRA image. Vessels in the neck and cerebral circulations are scaled from Melis et al. while vessels in the abdomen and lower body are scaled from Puelz et al. \cite{Puelz2017}. Figure \ref{fig:Network} shows the extracted geometries and extrapolated network used in this study. As suggested by Pennati et al. \cite{Pennati2000}, we scale vessel length using the following formula:
\begin{equation}
    L_2 = L_1\left(\frac{W_1}{W_2}\right)^\alpha,
    \label{eq:scale}
\end{equation}
where $\alpha = 0.35$ and $W_1$ (kg) and $L_1$ (cm) are patient weight and vessel length obtained from literature \cite{Melis2019,Puelz2017} respectively, and $W_2$ (kg) is the weight of the DORV or HLHS patient. The quantity $L_2$ (cm) is the  unknown vessel length. The inlet radii of vessels immediately outside of the imaged region are matched to the size of the outlet radii of their preceding vessels within the imaged region. Radii of downstream vessels are scaled by the same method as described above.
 \begin{figure}[t]
    \centering
    \includegraphics[width=0.75\textwidth]{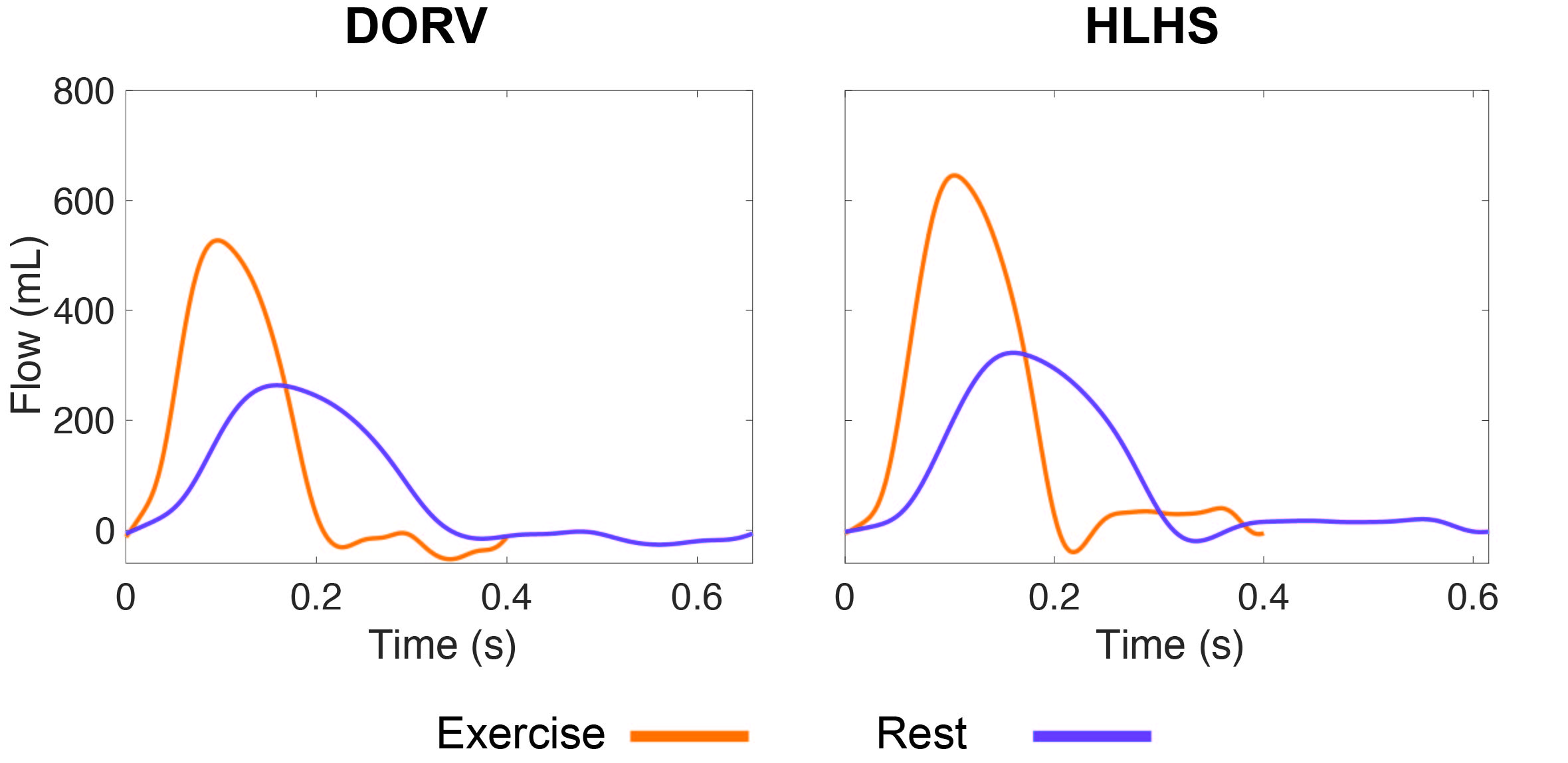}
    \caption{Inflow profiles for the DORV and HLHS patients at rest and light upright exercise conditions. Inflows at rest are extracted from 4D-MRI data. To simulate light upright exercise, inflow waveforms are multiplied by two and the duration of the cardiac cycle is shortened. At rest, DORV flow reaches about 300 mL at the inlet and under exercise conditions the inflow is doubled to reach about 600 mL at the inlet. At rest, HLHS flow reaches about 390 mL and under exercise conditions is doubled to about 680 mL.}
    \label{fig:Inflow}
\end{figure}

Large systemic vessels are known to taper along their length \cite{Caro2012}. To account for tapering in our models, we estimate tapering parameters by fitting an exponential function of the following form to data from the descending thoracic aorta:
\begin{equation}
\label{eq:taper}
    r = n_1\exp(-n_2x)+n_3.
\end{equation}
 In equation \eqref{eq:taper},  $r$ denotes the vessel radius (cm), $n_1$ denotes the inlet radius minus the outlet radius, $n_3$ denotes the outlet radius, $x$ is the axial location along the vessel, and $n_2$ denotes the average degree of taper. The estimated parameters are used to determine taper in vessels outside the imaged region.  Table \ref{table:Dim} lists the network and vessel dimensions. Inlet flow waveforms at rest (described in detail in \ref{Sec:Data}) for each patient were extracted from the 4D-MRI data.  Figure \ref{fig:Inflow} depicts the inlet flow waveforms for each patient at rest (measured) and at light upright exercise (estimated). The estimation of the inflow waveform during light upright exercise is done as follows. The waveform at the root of the ascending aorta was multiplied by two and the duration of the cardiac cycle was reduced by approximately 40\%, resulting in an average heart rate ($\sim$150 bpm) during exercise for females aged 8-18 years \cite{Harkel2011}.

\subsection{Fluid Dynamics Model}

Hemodynamics are predicted by solving a one-dimensional (1D) fluid dynamics model that describes blood flow, blood pressure, and vessel cross-sectional area along the vessel's axial dimension in a network of compliant vessels. Solutions to this model are computed explicitly in large vessels within the network, shown in Figure \ref{fig:Network}. The microvasculature hemodynamics are predicted using a linearized model that is solved semi-analytically within a structured tree framework. For this study, the microvascular networks provide impedances that are used as boundary conditions for the large vessel network. Below, we describe the large and small vessel models. Parameters for each patient are given in Tables \ref{table:Dim} and \ref{tab:param}.

\subsubsection{Large Vessels}

The 1D model is derived from the Navier-Stokes equations, assuming axially symmetric Newtonian flow, cylindrical vessels, and that blood is incompressible, viscous, and homogeneous with constant density $\rho=1.057$ (g/cm$^3$) and viscosity $\mu=0.032$ (g/cm/s) \cite{Olufsen2000}. The flow $q(x,t)$ (mL/s), the pressure $p(x,t)$ (g/cm/s$^2$), and the cross-sectional area $A(x,t)$ (cm$^2$) in each vessel satisfies the following 1D mass conservation and momentum balance equations:
\begin{eqnarray}
   && \frac{\partial A}{\partial t}+\frac{\partial q}{\partial x} = 0 \label{eq:cons1}\\
&& \frac{\partial q}{\partial t}+\frac{\partial}{\partial x}\left(\frac{q^2}{A}\right)+\frac{A}{\rho}\frac{\partial p}{\partial x} = -\frac{2\pi \nu R}{\delta}\frac{q}{A}+gA\cos(\theta) \label{eq:moment1},
\end{eqnarray}
where $0\leq x \leq L$ is the axial position within the vessel. Here,  $\nu = \frac{\mu}{\rho}$ (cm$^2$/s) is the kinematic viscosity, $g$ (cm/s$^2$) is the gravitational acceleration, and $\theta$ is the angle between the vessel and the gravitational field. 

These equations are obtained using an averaged velocity profile $u_x(r,x,t) = q/A$  along each vessel and by imposing a Stokes boundary layer in equation \eqref{eq:moment1} of the form
\begin{equation}
\label{eq:stokes}
u_x(r,x,t) = \left\{
  \begin{array}{lr} 
      \displaystyle \bar{u}_x, \hspace{1cm} & r<R-\delta, \\
       \displaystyle\bar{u}_x\frac{(R-r)}{\delta}, \hspace{1cm} & R-\delta<r\leq R, 
      \end{array}
\right.
\end{equation}
where $R$ is the radius of the vessel, $\delta= \sqrt{\nu T/2\pi}$ (cm) is the boundary layer thickness, $T$ (s) is the cardiac cycle duration, and $u_x$ is the average velocity in the axial direction. The system of equations is closed via a pressure-area relationship that describes the vessel wall as a linear elastic membrane:
\begin{equation}
    p(x,t)-p_0=\frac{4}{3}\frac{Eh}{r_0}\left(1-\sqrt{\frac{A_0}{A}}\right), \label{eq:wallmodel}
\end{equation}
where $E$ (g/cm/s$^2$) is Young's modulus, $h$ (cm) is vessel wall thickness, $p_0$ (g/cm/s$^2$) is the reference pressure, and $A_0$ (cm$^2$) is cross-sectional area when the pressure equals its reference value. To account for radius dependent stiffening, we define:  
\begin{equation}
\label{eq:elastic}
    \frac{Eh}{r_0}=k_1\exp(-k_2r_0)+k_3,
\end{equation}
where $k_1$ (g/cm/s$^2$), $k_2$ (1/cm) and $k_3$ (g/cm/s$^2$) are constants. Simulated pressures are converted to mmHg using the conversion that 1 mmHg equals 1333.22 g/cm/s$^2$.

The model described above forms a hyperbolic system of PDEs. Therefore, boundary conditions are required at the inlet and outlet of each vessel. At the inlet of the arterial tree (ascending aorta), a flow waveform is imposed using the waveform extracted from the 4D-MRI data (refer to Figure \ref{fig:Inflow}). At the junctions, we enforce mass conservation and pressure continuity:
\begin{equation}
\label{eq:flow}
    q_p = q_{d1}+q_{d2},
\end{equation}
\begin{equation}
\label{eq:pressure}
    p_p = p_{d1} = p_{d2},
\end{equation}
where subscript $p$ refers to the parent vessel and $di, i=1,2$ refer to the two daughter vessels. 

In equation \eqref{eq:moment1}, the term $gA\cos(\theta)$ accounts for gravitational effects. When modeling patients in a supine position, $\theta = \frac{\pi}{2}$. All vessels are approximately at the same height as the heart, so gravity does not need to be considered. When modeling patients in an upright position, all vessels carrying blood upwards from the heart are assigned $\theta = -\pi$, since blood is working against gravity. All vessels carrying blood downwards from the heart are assigned $\theta = 0$, since blood moves with gravity in the positive direction. All vessels moving approximately parallel to the ground are assigned $\theta = \frac{\pi}{2}$.

Similar to our previous studies \cite{Olufsen2000,Olufsen2001,Chambers2020,Colebank2019,Colebank2021}, model equations are non-dimensionalized and solved using the two-step Lax-Wendroff method.

\subsubsection{Small Vessels}
\begin{figure}[t!]
    \centering
    \includegraphics[width=0.75\textwidth]{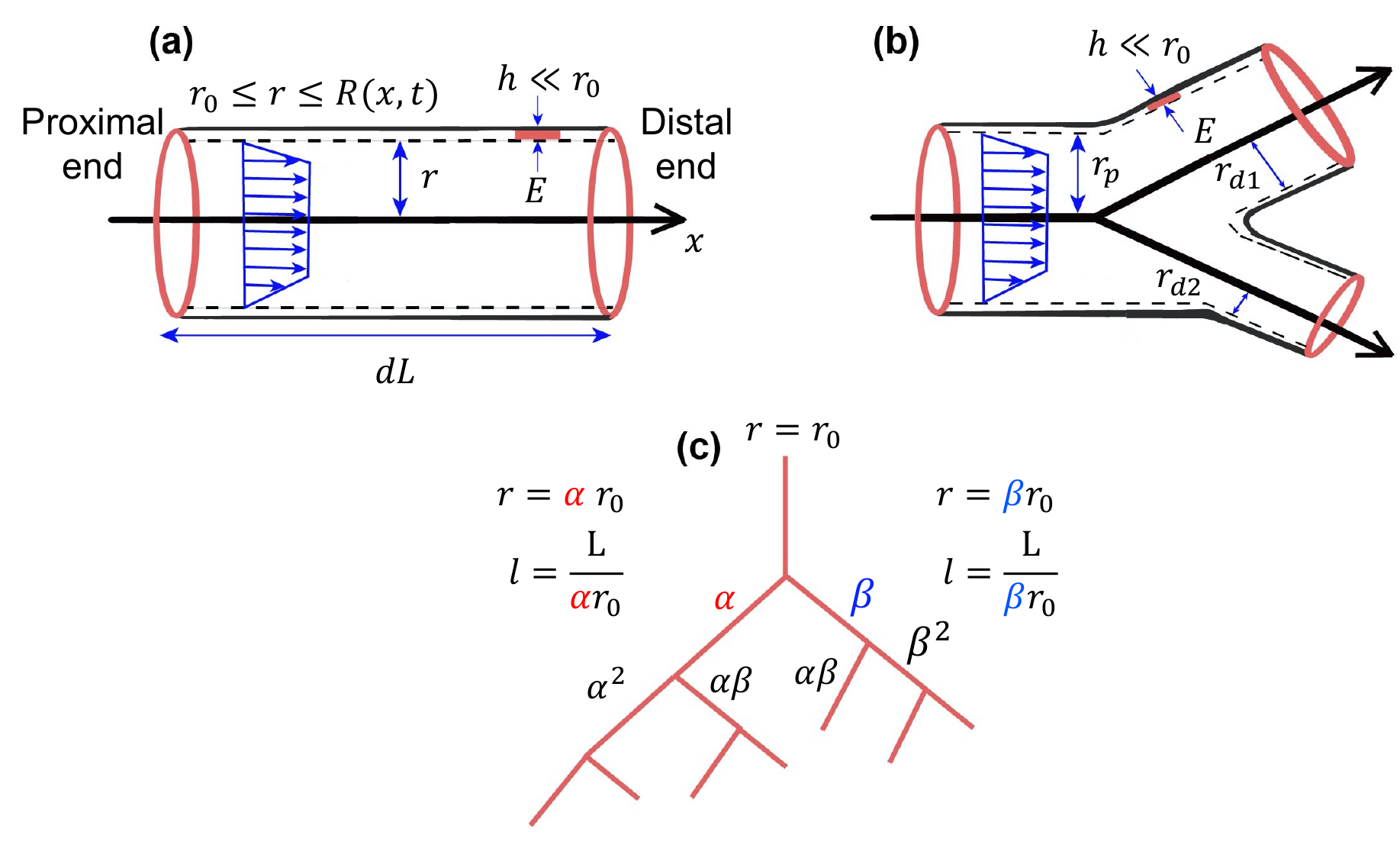}
    \caption{Illustration of the vessels modeled in this study. \textbf{(a)} shows the fluid dynamics of the large vessels. \textbf{(b)} shows the fluid dynamics of the large vessels at a junction to ensure continuity of pressure and conservation of flow. \textbf{(c)} shows the structured tree used the terminal vessels.}
    \label{fig:Model}
\end{figure}
Outflow boundary conditions are attached to each terminal vessel. Similar to our previous studies \cite{Chambers2020,Colebank2021,Olufsen2000,Olufsen2001}, we predict flow to the vascular beds by coupling large vessels to asymmetrical structured trees \cite{Olufsen2000,Olufsen2001}, as shown in Figure \ref{fig:Model}. In the small arteries, viscous forces are assumed dominate inertial forces, allowing us to linearize \eqref{eq:cons1} and \eqref{eq:moment1}. The equations are reduced by assuming periodicity of solutions, as described in detail in \cite{Olufsen2000}, resulting in:
\begin{eqnarray}
    &&i\omega Q + \frac{A_0(1-F_J)}{\rho}\frac{\partial P}{\partial x}=0 \label{eq:moment2},\\
    &&    i\omega CP + \frac{\partial Q}{\partial x} = 0,
    \label{eq:mass2}
\end{eqnarray}
where $F_j=2J_1(w_0)/w_0J_0(w_0), J_i(w_0)$, $i = 0,1$ are the zeroth and first order Bessel functions, and $w^2_0=i^3r^2_0\omega/v$. Note that $w^2=r^2_0\omega/v$ is the Womersley number. A reduced wave equation may be derived by differentiating equation \eqref{eq:mass2} and substituting the result into equation \eqref{eq:moment2}:
\begin{equation}
\label{eq:reducedwave}
    \frac{\omega^2}{c^2}Q+\frac{\partial^2Q}{\partial x^2}=0, \hspace{3mm} \text{or} \hspace{3mm} \frac{\omega^2}{c^2}P+\frac{\partial^2P}{\partial x^2} =0,
\end{equation}
where $c = \sqrt{A_0(1-F_j)/\rho C}$ is the wave propagation velocity. Solutions to equations \eqref{eq:moment2} and \eqref{eq:mass2} are:
\begin{equation}
\label{eq:soln1}
    Q = a\cos(\omega x/c)+b\sin(\omega x/c),
\end{equation}
\begin{equation}
\label{eq:sol2}
    P(x,\omega) = i\sqrt{\frac{\rho}{CA_0(1-F_j)}}(-a\sin(\omega x/c)+b\cos(\omega x/c)),
\end{equation}
where $a,b$ are constants of integration. As was done in the large vessels, boundary conditions need to be defined for the small vessels. Similarly, pressure continuity and mass conservation are preserved. A bifurcation is analogous to a transmission-line network in which the impedances satisfy:
\begin{equation}
\label{eq:impedance}
    \frac{1}{Z_p}=\frac{1}{Z_{d1}}+\frac{1}{Z_{d2}}.
\end{equation}
Similar to the original work by Olufsen et al. \cite{Olufsen2000}, we assume the terminal impedance at the end of the structured tree is zero.

\subsection{Wave intensity analysis}

To quantify the incident (forward-moving) and reflected (backward moving) components of waves, we use wave intensity analysis (WIA) \cite{Colebank2021,Qureshi2019,Qureshi2019a}. Assuming negligible frictional losses, the incident and reflected waves are approximated by setting $q=Au$ , where $u$ (cm/s) is fluid velocity, and defining:
\begin{equation}
\label{eq:wia1}
    \Gamma_{\pm} (t) = \Gamma_0+\int_0^Td\Gamma_{\pm}, \hspace{5mm} \Gamma=p,u
\end{equation}
\begin{equation}
\label{eq:wia2}
    dp_{\pm} = \frac{1}{2}(dp\pm\rho cdu), \hspace{3mm} du_{\pm}=\frac{1}{2}\left(du\pm\frac{dp}{\rho c}\right),
\end{equation}
where $c$ (cm/s) is the pulse wave velocity. We defined the time-normalized wave intensity as
\begin{equation}
\label{eq:wia3}
    WI_{\pm} = (dp_\pm/dt)(du_{\pm}/dt).
\end{equation}
Incident waves are classified as compressive with $WI_+, dp_+>0$ and expansive with $WI_+, dp_+<0$. Similarly, reflected waves are classified as compressive with $WI_-, dp_->0$ and expansive with $WI_-, dp_-<0$. The wave reflection coefficient is defined as the ratio of amplitudes of the reflected compression pressure waves to the incident compression pressure waves:
\begin{equation}
\label{eq:wia4}
I_R = \frac{\Delta p_-}{\Delta p_+}.
\end{equation}
Forward compression waves (FCW) begin at the inlet of the vessel and increase the pressure and flow velocity as blood propagates towards the outlet of the vessel. Forward expansion waves (FEW) also begin at the inlet of the vessel but decrease pressure and flow velocity. Backward compression waves (BCW) begin at the outlet of the vessel and propagate backward down the vessel, increasing pressure but decreasing flow velocity. Backward expansion waves (BEW) decrease pressure while accelerating flow \cite{Colebank2021}.

\subsection{Wall shear stress}
\label{sec:5}
We compute wall shear stress (WSS) in the large vessels, i.e., the stress the fluid exerts on the vessel wall, denoted $\tau_w$ (g/cm/s$^2$), by using the Stokes boundary layer given in equation \eqref{eq:stokes} \cite{Bartolo2022}. This results in the following equations: 
\begin{equation}
    \tau_w = -\mu \frac{\partial u}{\partial r},
\end{equation}
\begin{equation}
   \implies \tau_w = \left\{
  \begin{array}{lr} 
      0, \hspace{1cm} & r<R-\delta, \\
      \displaystyle\frac{\mu\bar{u}}{\delta}, \hspace{1cm} & R-\delta<r\leq R, 
      \end{array}
\right.
\end{equation}
where $\mu$ is blood viscosity and $\delta = \sqrt{\nu T/2\pi} $ is the boundary layer thickness.

\subsection{Model calibration and simulation}
\label{sec:6}
Dimensions from Table \ref{table:Dim} and inlet flow waveforms (cardiac output and heart rate) mentioned previously (Figure \ref{fig:flowdata}) are used to calibrate models for each patient. The same connectivity matrix is defined for both patients, but the vessel length and radii are different.  We assume that neither patient has microvascular disease at the time of the measurements. Therefore, both patients are assigned the same boundary condition parameters, i.e. the structured tree parameters including the fractal scaling ratios $\alpha$ and $\beta$, the minimum radius $r_\text{min}$, and the length to radius ratio lrr.  Nominal large vessel stiffness parameters from Olufsen et al. \cite{Olufsen2000} are tuned for each patient to generate flow waveforms and pressure predictions that match both the data from the imaged region and the cuff pressures. The tuned values are listed in Table \ref{tab:param}. Small vessel stiffness, denoted $ks_i$ in Olufsen et al. \cite{Olufsen2000}, are equal to their large vessel counterparts (i.e., $k_1 = ks_1$, etc.).

The calibrated model is used to predict $p(x,t)$, $q(x,t)$, and $A(x,t)$ in all 57 vessels in the supine position and during light upright exercise. Using these predictions, we perform WIA to compare forward and backward wave propagation in the two patients. From the flow and area predictions, we calculate WSS in the large vessels. WSS sensed by endothelial cells in the arteries is typically altered in cardiovascular disease, making it an important output quantity to consider.

\begin{table}[b!]
\caption{Stiffness and structured tree parameters for the arterial network. Parameter values were nominally set to be consistent between DORV and HLHS patients. Specific values were found by manual tuning for particular vessels to fit flow predictions to patient data.}
\label{tab:param}
\centering
\adjustbox{max width = \textwidth}{
\begin{tabular}{lrr}
\hline\noalign{\smallskip}
\multicolumn{3}{c}{\textbf{Nominal Parameter Values}}\\
\noalign{\smallskip}\hline\noalign{\smallskip}
 $k_1$  & 2.0 $\times$ 10$^6$ & g/s$^2$cm\\
 $k_2$  & -35 & cm$^{-1}$\\
 $k_3$  & 3.80 $\times$ 10$^5$ & g/s$^2$cm\\
 $r_\text{min}$ & 0.01 & cm \\
 $\alpha$  & 0.90 &\\ 
 $\beta$  & 0.60 &\\
 lrr  & 50 &\\
\noalign{\smallskip}\hline\noalign{\smallskip}
\multicolumn{3}{c}{\textbf{Tuned} $r_\text{min}$ \textbf{(cm)}}\\
\noalign{\smallskip}\hline\noalign{\smallskip}
\textbf{Vessel} & \textbf{DORV} & \textbf{HLHS}\\
 10,17 & 0.03 & 0.03\\
 12,14 & 0.001 & 0.001\\
 21,23 & 0.001 & 0.001\\
 27,29 & 0.001 & 0.001\\
 31,33 & 0.001 & 0.001\\
 46,48 & 0.01 & 0.10\\
 47,49 & 0.01 & 0.10\\
 50,52 & 0.01 & 0.10\\
 51,53 & 0.01 & 0.10\\
\noalign{\smallskip}\hline\noalign{\smallskip}
\multicolumn{3}{c}{\textbf{Tuned $k_3$ (g/cm/s$^2$)}}\\
\noalign{\smallskip}\hline\noalign{\smallskip}
\textbf{Vessel} & \textbf{DORV} & \textbf{HLHS}\\
 1 & 3.8 $\times$ 10 $^5$ & 5.7 $\times$ 10 $^5$ \\
 2 & 3.8 $\times$ 10 $^5$ & 5.7 $\times$ 10 $^5$ \\
 4 & 3.8 $\times$ 10 $^5$ & 5.7 $\times$ 10 $^5$ \\
 5,16 & 3.8 $\times$ 10 $^5$ & 4.56 $\times$ 10 $^5$\\
 6,11 & 7.6 $\times$ 10 $^5$ & 1.9 $\times$ 10 $^6$\\
 7,9 & 3.8 $\times$ 10 $^5$ & 5.7 $\times$ 10 $^5$\\
 8 & 3.8 $\times$ 10 $^5$ & 5.7 $\times$ 10 $^5$ \\
 10,17 & 7.6 $\times$ 10 $^5$ & 3.8 $\times$ 10 $^5$\\
 12,14 & 2.66 $\times$ 10 $^6$ & 2.66 $\times$ 10 $^6$\\
 13,15 & 2.66 $\times$ 10 $^6$ & 2.66 $\times$ 10 $^6$\\
 21,23 & 1.9 $\times$ 10 $^6$ & 2.66 $\times$ 10 $^6$\\
 35,37,39,41,43 & 3.8 $\times$ 10 $^5$ & 3.04 $\times$ 10 $^5$\\
 44-52 & 3.8 $\times$ 10 $^5$ & 3.04 $\times$ 10 $^5$\\
\noalign{\smallskip}\hline
\end{tabular}}
\end{table}

\section{Results}

Results are computed for the DORV and HLHS patients using parameters listed in Table \ref{table:Patient}. Pressure and flow predictions are shown in the aorta and regions of interest, including the brain and gut (Figures \ref{fig:Local} and \ref{fig:Total}). WIA and WSS results are shown along the aorta (Figures \ref{fig:WIA} and \ref{fig:WSS}). 

\label{sec:7}
\subsection{Pressure and flow predictions}
\begin{figure}[b!]
    \centering
    \includegraphics[width=1.0\textwidth]{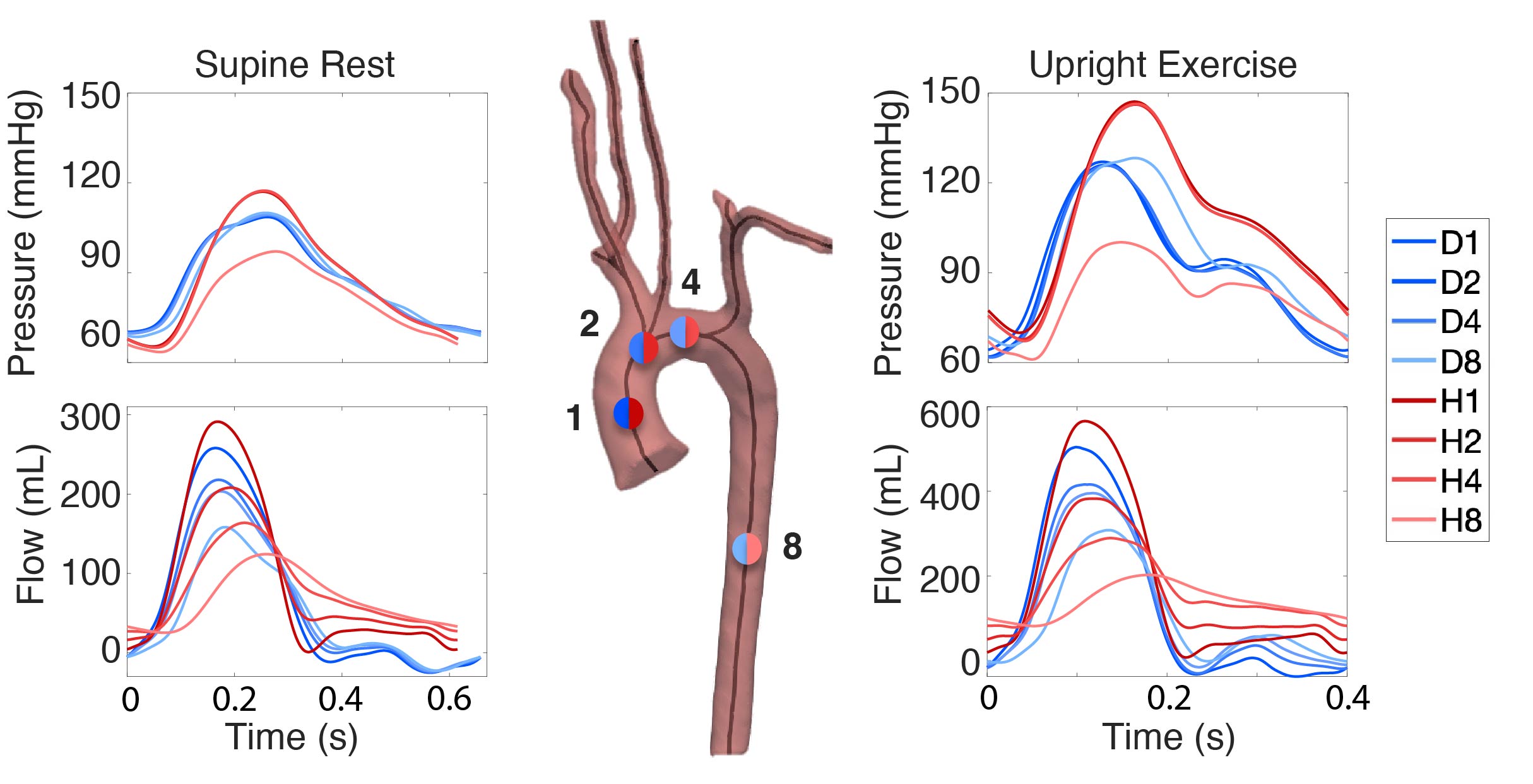}
    \caption{Localized pressure and flow predictions in the imaged region. Pressures for the DORV patient range from ~65-105 mmHg at supine rest. A greater range is seen in the HLHS patient, from ~63-118mmHg at supine rest. Cardiac output is similar in the aorta in supine rest and light upright exercise conditions. Predictions deviate significantly in the distal part of the descending aorta. D1, D2, D4, and D8 denote predictions for the DORV patient. H1, H2, H4, and H8 denote predictions for the HLHS patient.}
    \label{fig:Local}
\end{figure}

\begin{figure}[t]
    \centering
    \includegraphics[width=1.0\textwidth]{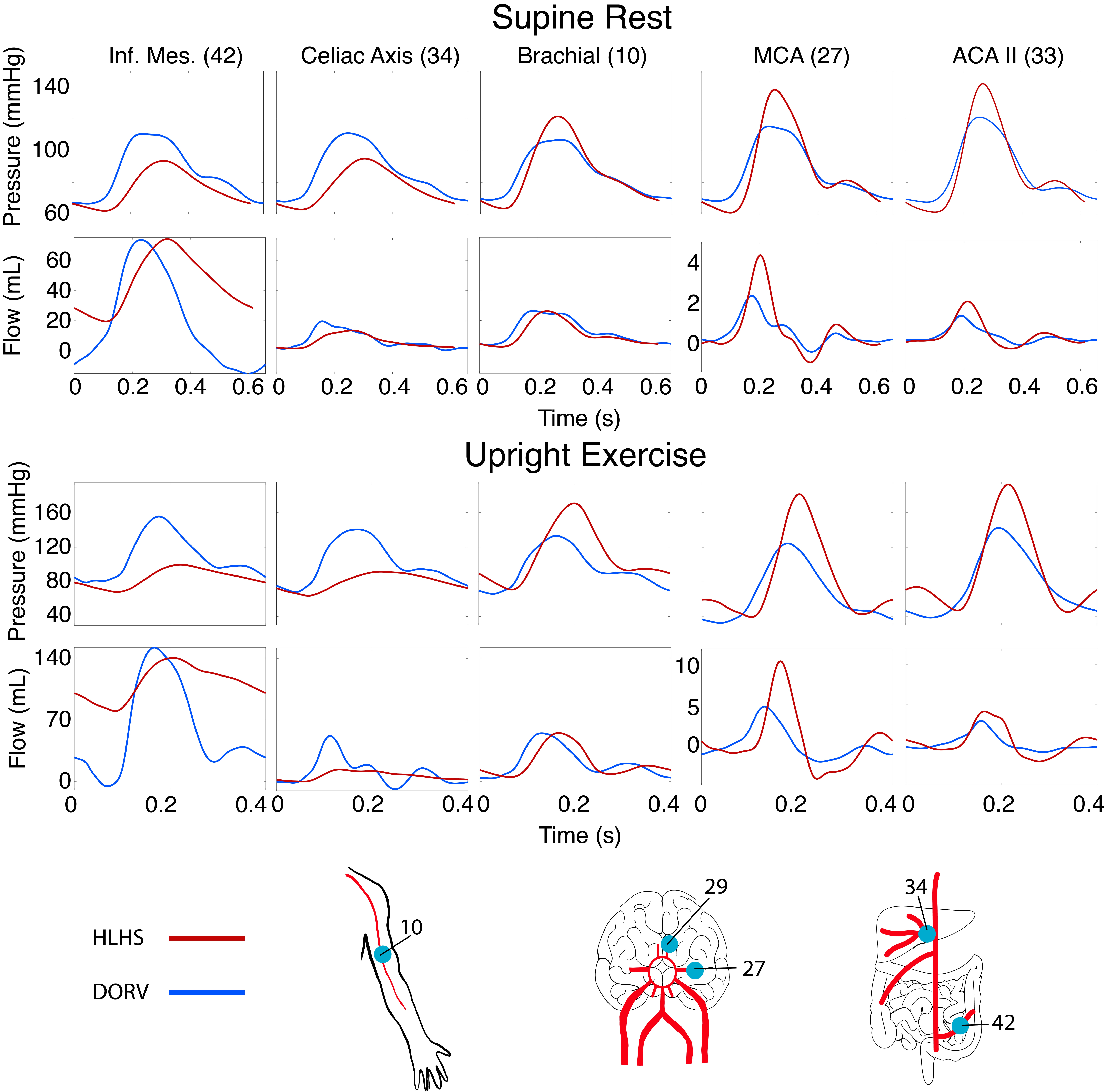}
    \caption{Pressure and flow predictions in the supine rest and light upright exercise positions for different vessels. Cerebral pressures are greater in the HLHS patient but flow is similar to the DORV patient at rest and exercise. Pressure and flow in the gut are reduced in the HLHS patient compared to the DORV patient.}
    \label{fig:Total}
\end{figure}

Simulations under supine rest and light upright exercise conditions are shown in Figures \ref{fig:Local} and \ref{fig:Total}. Figure \ref{fig:Local} shows predictions within the segmented network and Figure \ref{fig:Total} shows predictions in the peripheral network which describes systemic arterial vessels outside of the imaged region. 
\begin{table}[b]
\caption{Fraction of blood volume flowing to the brain, liver and gut, and lower body.}
\centering
\label{tab:flow}   
\begin{tabular}{l|ll|ll}
\hline\noalign{\smallskip}
& \multicolumn{2}{|c|}{\textbf{Rest}} & \multicolumn{2}{|c}{\textbf{Exercise}}\\
\noalign{\smallskip}\hline\noalign{\smallskip}
 \textbf{Region} & \textbf{DORV} & \textbf{HLHS} & \textbf{DORV} & \textbf{HLHS}\\
\noalign{\smallskip}\hline\noalign{\smallskip}
Cerebral & 0.086 & 0.073 & 0.048 & 0.052\\
Liver and Gut & 0.243 & 0.157 & 0.202 & 0.094\\
Lower Body & 0.257 & 0.480 & 0.400 & 0.608\\
\noalign{\smallskip}\hline
\end{tabular}
\end{table}

Model predictions of the systolic/diastolic pressures for the DORV and HLHS patients are 109/68 mmHg and 118/63 mmHg respectively. Pressures for both patients are in the normotensive range, but the HLHS patient has significantly higher pulse pressure. These pressure predictions agree with measured values listed in Table \ref{table:Patient}, which are taken in the supine position. The cardiac output, which is 5.08 L/min for the DORV patient and 4.06 L/min for the HLHS patient, is also within the normal range for a healthy bi-ventricular heart. Even though cardiac output is lower for the HLHS patient, as shown in Figure~\ref{fig:Inflow}, the pulse flow is similar for the two patients. Since this study imposes flow at the inlet of the ascending aorta, cardiac output is a model input. The same applies to heart rate, which is 91 bpm and 97 bpm for the DORV and HLHS patients respectively. Predicted blood pressure and flow waveforms in the aorta are shown in Figure~\ref{fig:flowdata}, demonstrating that the calibrated model fits measurements for both patients. An exception is in the thoracic aorta for the HLHS patient, where the model predicted pulse flow is less than that of the data, but the average flow from the model and data agree relatively well.

In light upright exercise, differences in blood flow through the aorta become more apparent. Even though the aortic root pulse flow is similar between patients, the descending aortic pulse flow in the HLHS patient is significantly lower compared to the DORV patient (Table \ref{tab:flow}). Figure \ref{fig:Total} demonstrates that for both patients, cerebral flow in the MCA and ACA II is similar in supine rest and light upright exercise. However, there is less flow to the lower extremities in the HLHS patient at rest compared to the DORV patient. The difference between the DORV and HLHS patient is more dramatic during light upright exercise for flow to the gut and lower extremities. Table \ref{tab:flow} lists the fraction of blood volume entering these vessels. Values were calculated by dividing the average blood flow through each vessel by the average blood flow through the ascending aorta. 

\subsection{Wave intensity analysis}
\label{sec:9}
    \begin{figure}[b!]
        \centering
        \includegraphics[width=1.0\textwidth]{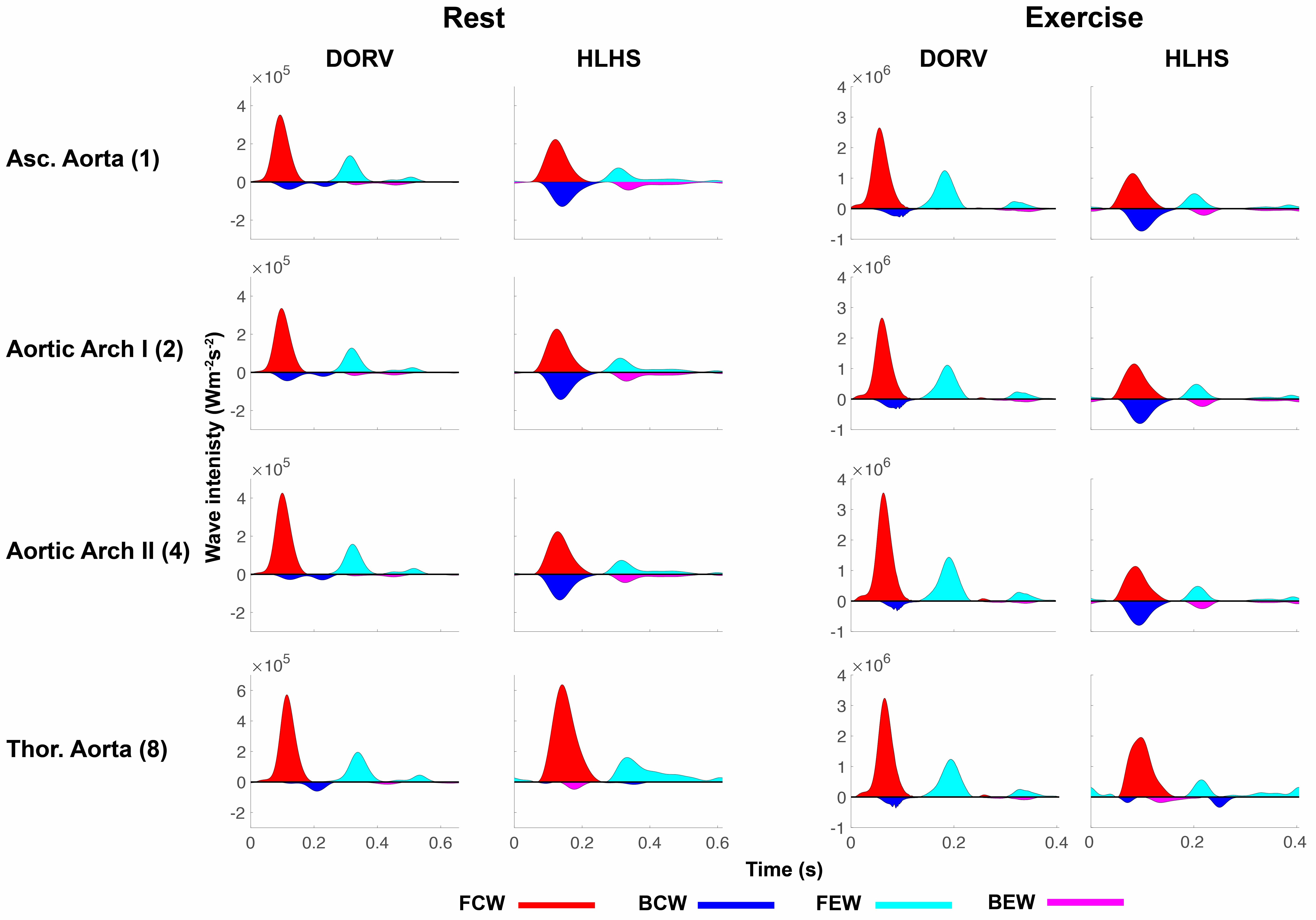}
        \caption{Wave intensity analysis at the midpoint of each aortic vessel segment at supine rest and light upright exercise. Results for both patients show a predominant forward traveling wave in each vessel. The HLHS patient has larger backward traveling waves while the DORV patient has relatively minor backward traveling waves.}
        \label{fig:WIA}
    \end{figure}
    
Wave intensity profiles for both patients at the midpoint of the aorta are shown in Figure \ref{fig:WIA} for both supine rest and light upright exercise conditions. In comparison to the other waves present, each vessel in all cases has a predominant FCW propelling blood down the vessel. In the DORV patient at rest and light upright exercise, the BCWs are small in the aorta. Exercise predictions for the HLHS patient differ from the predictions for the DORV patient in that the BCW is larger, except for wave in the thoracic aorta. The FCW is also decreased compared to the DORV patient. The DORV patient has an increased FEW compared to the HLHS patient at both rest and light upright exercise. The wave reflection coefficients are shown in Table \ref{tab:WaveReflec}. The DORV patient has a consistently smaller wave reflection coefficient compared to the HLHS patient at both rest and light upright exercise.

\begin{table}[t!]
\caption{Wave reflection coefficients ($I_R$) computed for the aortic vessels.}
\centering
\label{tab:WaveReflec}   
\begin{tabular}{l|ll|ll}
\hline\noalign{\smallskip}
& \multicolumn{2}{|c|}{\textbf{Rest}} & \multicolumn{2}{|c}{\textbf{Exercise}}\\
\noalign{\smallskip}\hline\noalign{\smallskip}
 \textbf{Region} & \textbf{DORV} & \textbf{HLHS} & \textbf{DORV} & \textbf{HLHS}\\
\noalign{\smallskip}\hline\noalign{\smallskip}
Ascending aorta & 0.567 & 0.834 & 0.359 & 0.816 \\
Aortic arch I & 0.583 & 0.830 & 0.389 & 0.835 \\
Aortic arch II & 0.483 & 0.788 & 0.310 & 0.819 \\
Thoracic aorta & 0.389 & 0.170 & 0.377 & 0.347 \\
\noalign{\smallskip}\hline
\end{tabular}
\end{table}

\subsection{Wall shear stress}
\label{sec:10}
\begin{figure}[b!]
    \centering
    \includegraphics[width=0.75\textwidth]{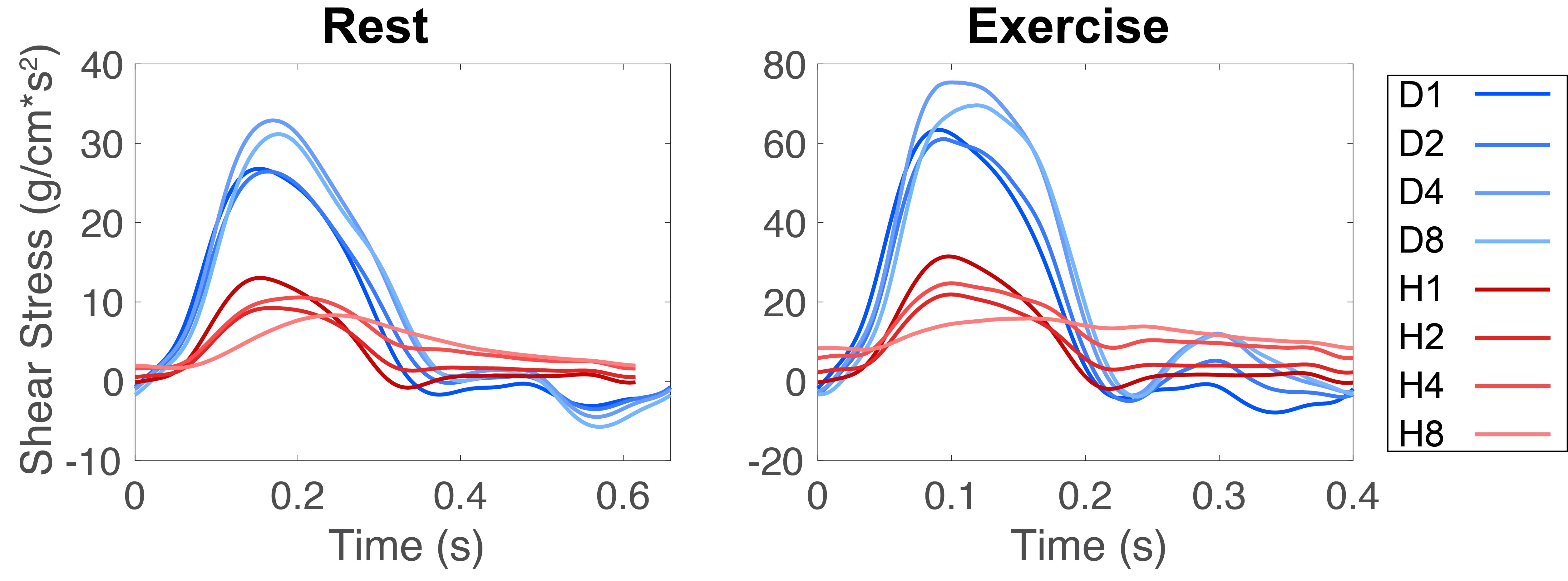}
    \caption{WSS results within the aortic vessels. The DORV patient is shown in blue and the HLHS patient is shown in red. During supine rest, WSS values for the DORV patient have a maximum of $\sim$30-35 g/cm/s$^2$ while values for the HLHS patient have a maximum of $\sim$10-15 g/cm/s$^2$. During light upright exercise, WSS values for the DORV patient have a maximum of $\sim$75-80 g/cm/s$^2$ while WSS values for the HLHS patient have a maximum of $\sim$30-40 g/cm/s$^2$.}
    \label{fig:WSS}
\end{figure}

WSS results in the aorta are shown in Figure \ref{fig:WSS}. At rest, the DORV patient has a maximum shear stress of $\sim$30-35 g/cm/s$^2$ and the HLHS patient has a maximum shear stress of $\sim$10-15 g/cm/s$^2$. During light upright exercise, these maximums are essentially doubled; the DORV patient has a maximum of $\sim$75-80 g/cm/s$^2$ and the HLHS patient has a maximum of $\sim$30-40  g/cm/s$^2$.

\section{Discussion}
\label{sec:11}

In this study, we described the construction of models to predict blood pressure and flow dynamics in an HLHS patient with a reconstructed aorta and a DORV patient with a native aorta. Cardiac output was lower and pulse pressure higher in the HLHS patient as compared to the DORV control patient. This prediction was expected in part due to the remodeling and abnormal morphology of the reconstructed aorta, but data for more patients will be necessary in order to draw significant conclusions. Models were calibrated to flow and pressure data from both patients. Calibrated flow predictions agreed with measured waveforms extracted from 4D-MRI data (Figure \ref{fig:flowdata}). Pressure predictions were calibrated to sphygmomanometer measurements taken in the supine position. 



\subsection{Pressure and flow predictions}

Pressure and flow predictions were generated in vessel networks that were extracted from MRA data. Patient-specific inlet flow waveforms derived from 4D-MRI images were used as boundary conditions at the inlet of the ascending aorta (Figure \ref{fig:flowdata}). Parameters determining vessel stiffness and peripheral vascular resistance were tuned to fit measurements, and model outcomes were calculated using the calibrated models. 

Model predictions revealed that blood flow to the liver and gut circulations were significantly lower in the HLHS patient, while flow through the head and neck vessels were similar between patients.  The decrease in flow to the lower body was a result of less flow through the aorta (Figure \ref{fig:Local}). Decreased perfusion to the gut for the HLHS patient might explain the more frequent development of FALD in this population \cite{Camposilvan2008}. Decreased flow to the gut for the HLHS patient may also be consistent with the body's adaptative mechanisms that ensure the brain is supplied with an adequate amount of blood \cite{Saiki2016}. Previous studies, including by Navaratnam et al. \cite{Navaratnam2016}, have shown that HLHS patients with reconstructed aortas eventually have inadequate perfusion to the brain. In particular, the authors identified a group of adults with a Fontan circulation that experienced reduced oxygen delivery to the systemic circulation under exercise conditions. Furthermore, they found higher cerebral deoxygenation during vigorous exercise. Our model did not show a significant difference in perfusion to the brain between the two patients. This discrepancy might result from our choice to only simulate light upright exercise and to not include the impact of changes in peripheral vascular resistance, compliance, and cardiac contractility imposed by the autonomic control system during vigorous exercise. 

Blood pressure for the HLHS patient was higher in the head and neck vessels, compared to the gut vessels. These results suggest hypertension according to a study by Blanco et al. \cite{Blanco2017} that predicted cerebral pressures in both normotensive and hypertensive double-ventricle patients. Their normotensive studies predicted a peak pressure of about 100 mmHg in the middle cerebral artery (MCA). Our pressure predictions for the DORV patient in the MCA and ACA II were consistent with these results; our model predicted a peak systolic pressure of approximately 100-115 mmHg. At rest, the systolic pressure predictions for the HLHS patient were approximately 140 mmHg in both the MCA and anterior cerebral artery (ACA) II, with the diastolic pressures around 58 mmHg. These values are slightly hypertensive and could be indicative of peripheral artery stiffening \cite{Qamari2020}. This result was exacerbated under exercise conditions; the systolic blood pressure in the MCA and ACA II for the HLHS patient was ~190 mmHg. A review by Qamari et al. \cite{Qamari2020} found that increased blood pressure not only can lead to remodeling as a compensatory mechanism but can also lead to ischemic damage to the brain and an increased risk of stroke. It is not uncommon for HLHS patients to experience cerebral circulation issues. A study by Kotani et al. \cite{Kotan2018} investigated causes of death from a subgroup of Fontan patients over 20 years. Deaths due to cerebral issues, including thromboembolism, were common in patients who survived the Fontan operation. Similar 1D modeling studies, including the work by Puelz et al. \cite{Puelz2017}, demonstrated similar flow and pressure predictions in a Fontan circulation for an HLHS patient. Their results included a maximum flow in the superior mesenteric artery of approximately 15 mL/s. Peak systolic flow in the superior mesenteric artery of our model was around 13 mL/s (data not shown).


\subsection{Wave reflections and wall shear stress}

WIA is an important tool for assessing circulatory function \cite{Broyd2015}. Results from the HLHS patient in both rest and light upright exercise showed smaller FCWs and greater BCWs. Exercise conditions made these differences more apparent. Overall, backward traveling waves, which decelerate the flow (see e.g. \cite{Broyd2015}), were of higher magnitude in the HLHS patient. Pomella et al. \cite{Pomella2018} compared WIA results from healthy individuals during rest and exercise. Their results showed an increase in forward traveling waves and a decrease in backward traveling waves during exercise. These findings are consistent with our results for the DORV patient, suggesting that the HLHS patient may have abnormal circulatory function during light upright exercise. Our results for the HLHS patient appear consistent with the results from Schafer et al. \cite{Schafer2021}. Their study found that HLHS reconstructed patients have decreased FCW and an increased BCW/FCW ratio compared to Fontan patients with other single left ventricle diseases. Values of the wave reflection coefficient, ($I_R$), provide a measure of local wave reflections between the vessel and the downstream vasculature. Our results showed that upright exercise consistently reduced the reflection coefficient in nearly all vessels for both patients, except for the thoracic aorta in the HLHS patient. Though we do not account for the acute regulatory effects of exercise in our model (e.g., peripheral vasodilation), we still observed a decrease in the wave reflection coefficient that is consistent with prior studies \cite{Pomella2018}. 

WSS can provide insight into endothelial mechanotransduction and long-term adaptation due to hyper/hypotension. WSS results for the DORV patient at rest were typical of a healthy individual, with a maximum of 30-35 \cite{Callaghan2018}. The maximum WSS for the HLHS patient is reduced, which can be indicative of hypertension and vessel wall deformation \cite{Traub1998}. Yang et al. \cite{Yang2014} investigated the correlation between WSS and hypertension to ultimately determine if local WSS values can be associated with vascular deformation. They found that hypertensive patients had lower peak WSS values due to an increased arterial diameter as a compensatory mechanism to combat stiffening vessels. The increase in arterial diameter led to stagnate blood flow in patients in their study, as seen in HLHS patients post-Norwood procedure (due to an increased aortic diameter). Low WSS also caused the transcription of genes to downregulate nitric oxide and upregulate endothelin-1, causing vasoconstriction, increased blood pressure, and led to further degradation of the vessel wall by promoting smooth muscle cell growth and causing loss of vessel compliance \cite{Traub1998}.

Predictions of shear stress have become increasingly important for optimizing patient outcomes. A recent study by Loke et al. \cite{Loke2020} combined model simulations with expert surgeon input to develop a Fontan conduit that minimized hemodynamic forces. The study used computational modeling across several geometries to help minimize power loss, improve hepatic flow, and minimize WSS. However, few have computationally investigated arterial WSS in Fontan patients, as we have done here. Forecasts of WSS from computational models can assist in determining possible interactions between endothelial cell dysfunction and poor outcomes post Fontan surgery.

\subsection{Limitations}

The computational efficiency of the 1D model used in this paper allows for rapid calibration of the model to clinical data. Another advantage of our modeling approach is the ability to predict flow in the entire peripheral vessel network, including vessels outside of the imaged region. However, this model also has several limitations, particularly for HLHS patients in which the reconstructed aorta is not cylindrical. In reality, the abnormal morphology of the reconstructed aorta likely contributes to energy loss due to the formation of vortices.  A description of energy loss can be included in the 1D model, as was done in several previous studies \cite{Colebank2021,Mynard2015}. An energy loss model could then be calibrated by comparing 1D and 3D simulations, e.g. using a 3D fluid-structure interaction (FSI) model for the reconstructed aorta \cite{Baumler2020,Bazilevs2009,Griffith2020}, or by identifying potential regions with secondary flow directly from 4D-MRI data. Another limitation is the lack of a heart model in our study. We did not include a description of the single ventricle, and instead we used a prescribed inflow waveform derived from 4D-MRI data.  Also, we extended the vessel network beyond the imaged region by scaling vessel lengths and radii based on a healthy individual's geometry. However, DORV and HLHS patients could have abnormalities in their systemic circulation. Lastly, the scaling performed for the upright exercise model corresponded to an average healthy female adolescent. However, heart rate and cardiac output in exercise for HLHS and DORV patients might change to a different degree compared to healthy individuals.

Another limitation is the lack of quantification in the uncertainty and measurement error of our data. Flow velocities derived from the 4D-MRI data were averaged over several cardiac cycles. In addition, flow was not exactly conserved in the original data. This lack of conservation was particularly evident for the HLHS patient. Furthermore, flow within the the left common carotid of the DORV patient was lower than the typical physiologic range. These discrepancies might be attributed to separate masks that were used for the ascending/descending aorta and head and neck vessels in order to post-process the 4D-MRI data and generate the flow waveforms. Finally, we had a single pressure measurement from a blood pressure cuff. For better model calibration, it would be useful to have several pressure readings in different vessels or multiple cuff pressure measurements. 



\section{Conclusion}

This study described the construction of models that predict blood pressure and blood flow in a 57-vessel network for the systemic arteries. Models were constructed from MRA data and calibrated to 4D-MRI and sphygmomanometer pressure data for both a DORV patient and an HLHS patient. The calibrated models were used to compute flow and pressure waveforms, perform wave intensity analysis, and predict wall shear stress at rest and during light upright exercise. The two patients in our study had normal cardiac output, blood pressure and heart rate. The aorta for the HLHS patient was significantly wider than the control patient and the wall stiffness determined from model calibration was higher. Our results showed that the HLHS patient had decreased flow to the gut and increased cerebral pressures, and these findings were consistent with the literature. These results were more prominent during light upright exercise. Both patients had similar flow to the brain at rest and during light upright exercise. Results from wave intensity analysis suggested abnormal flow in the HLHS patient with similar magnitude forward and backward traveling waves. Wall shear stress results for the HLHS patient showed low values in the aortic vessels, suggesting the patient may have had hypertension.

\label{sec:12}
\begin{acknowledgements}
This material is based upon work supported by the National Science Foundation Graduate Research Fellowship under Grant No. DGE-2137100. Any opinion, findings, and conclusions or recommendations expressed in this material are those of the authors(s) and do not necessarily reflect the views of the National Science Foundation. The project described was supported by the National Center for Research Resources and the National Center for Advancing Translational Sciences, National Institutes of Health, through Grant TL1 TR001414 (MJC). The content is solely the responsibility of the authors and does not necessarily represent the official views of the NIH. 

\end{acknowledgements}



\end{document}